\documentclass[letterpaper,10pt,conference,draftcls,onecolumn]{IEEEtran}
\IEEEoverridecommandlockouts
\usepackage{cite}
\usepackage{amsmath,amssymb,amsfonts,amsthm}
\usepackage[ruled,norelsize]{algorithm2e}
\usepackage{graphicx}
\usepackage{textcomp}
\usepackage{xcolor}
\usepackage{enumerate}
\usepackage{mathtools}
\usepackage{etoolbox}
\usepackage{balance}

\providebool{Report}
\setbool{Report}{true}

\newcommand{\R}{\mathbb{R}}
\newcommand{\N}{\mathbb{N}}

\newtheorem{assumption}{Assumption}
\newtheorem{theorem}{Theorem}
\newtheorem{problem}{Problem}
\newtheorem{lemma}{Lemma}
\newtheorem{definition}{Definition}

\newtheorem{proposition}{Proposition}

\newcommand{\col}{\textnormal{col}}

\usepackage{mathtools}
\mathtoolsset{showonlyrefs=true}

\usepackage{color}

\usepackage{marginnote}

\def\BibTeX{{\rm B\kern-.05em{\sc i\kern-.025em b}\kern-.08em
    T\kern-.1667em\lower.7ex\hbox{E}\kern-.125emX}}
\begin{document}

\title{Exponentially Converging Distributed Gradient Descent with Intermittent Communication via Hybrid Methods\\

\author{
Katherine R. Hendrickson$^a$\thanks{${}^a$K. Hendrickson and M. Hale are with the Department of Mechanical and Aerospace Engineering, University of Florida. Emails:
\texttt{\{kat.hendrickson,matthewhale\}@ufl.edu}. 
Research partially supported by AFOSR under Grant no. FA9550-19-1-0169, by 
ONR under Grant no. N00014-19-1-2543, 
and by a task order contract from the Air Force Research Laboratory through Eglin AFB.
}, Dawn M. Hustig-Schultz$^b$\thanks{${}^b$D. Hustig-Schultz and R. G. Sanfelice are with the Department of Computer Engineering,
University of California, Santa Cruz. Emails:
\texttt{\{dhustigs,ricardo\}@ucsc.edu}. Research partially supported by NSF Grants no. ECS-1710621, CNS-1544396, and CNS-2039054, by AFOSR under Grants no. FA9550-19-1-0053, FA9550-19-1-0169, and FA9550-20-1-0238, by the Army Research Office under Grant no. W911NF-20-1-0253, and by CITRIS and the Banatao Institute at the University of California.
}, Matthew T. Hale$^a$, and Ricardo G. Sanfelice$^b$}
}

\maketitle

\begin{abstract}
We present a hybrid systems framework for multi-agent optimization in which
agents execute computations in continuous time and communicate in discrete
time. The optimization algorithm is a hybrid version of parallelized
coordinate descent. Agents implement a sample-and-hold strategy
in which gradients are computed at communication times and held
constant during flows between communications. 
Completeness of maximal solutions under these hybrid dynamics is established.
Under assumptions of smoothness and strong convexity, we show that this system exponentially converges to the minimizer of an objective function. 
Simulation results illustrate this convergence rate. 
\end{abstract}

\section{Introduction}
Convex optimization problems arise in many
areas of engineering, including machine learning~\cite{sra12},
communications~\cite{luo06}, robotics~\cite{verscheure09}, and others. 
Fundamentally, regardless of the application area, the goal is to design
an algorithm that will converge to a minimum of an objective function,
possibly under some constraints. Recently, 
there has been increased interest in studying 
optimization algorithms
in continuous time using
tools from dynamical systems to establish convergence
to minimizers; see~\cite{su16,rahili2017,garg20}.

In this paper, we develop a hybrid optimization algorithm for the analysis
of multi-agent systems with continuous-time
updates and intermittent discrete-time communication events. This is motivated
by two factors. First, we wish to leverage the large collection of 
tools from dynamical systems to analyze multi-agent optimization. 
Second, there exist many multi-agent controllers
that operate in continuous time to minimize some objective function,
e.g., in consensus~\cite{olfati07} and coverage control~\cite{cortes04},
and our analyses will apply to such systems. 
However, while
individual agents' computations occur in continuous time, communication between them inherently occur in discrete time because
communicated information arrives at isolated time instants. 
This mixture of continuous- and
discrete-time elements naturally leads us to a hybrid system model. 

The algorithm we propose is essentially a hybrid version 
of parallelized block coordinate descent~\cite{bertsekas89}, in which each agent
updates only a small subset of all decision variables 
in continuous time, and agents
communicate these updates to others in discrete time. In the proposed model, communication between agents
occurs when a decreasing timer reaches zero, at which point the timer
is reset to some value within a specified range. 
Agents use a sample-and-hold strategy in which gradients
are computed at the communication times and then held constant
and used continuously until the next communication event. 
This approach is inspired by recent work~\cite{PHILLIPS2019} that has
successfully applied it to synchronization problems. 
We consider objective functions that satisfy 
typical, mild assumptions for distributed optimization, namely
strongly convex objective functions with Lipschitz gradients. 

We leverage the theory of hybrid systems to prove that
the proposed hybrid algorithm has several desirable properties.
First, we define a hybrid system model for this algorithm
and show that, under these hybrid dynamics, every maximal
solution is complete, with domain allowing arbitrarily large ordinary time. As a result, there are no
theoretical obstructions to running this algorithm
for arbitrarily long periods of time.
Second, we use Lyapunov analysis to show that, even under intermittent information sharing,
the hybrid optimization algorithm exponentially converges to the minimizer of an objective function. Furthermore, we derive
an explicit convergence rate in terms of system parameters. 

The developments in this paper can be regarded as continuous-time counterparts to
``classical'' discrete-time algorithms in multi-agent 
optimization~\cite{bertsekas89}. 
Related research in multi-agent continuous-time optimization includes~\cite{gharesifard2014,lu2012,rahili2017}, though 
those works all use a consensus-based update law that executes computations
and communications both in continuous time. However, we avoid continuous-time
communications to account for cases in which they are not possible
or simply undesirable,
e.g., over long distances or when power is limited. 

The most similar works are~\cite{PHILLIPS2019},~\cite{KIA2015}, which also study continuous-time optimization
with discrete-time communication. However, those works also use consensus-based optimization
algorithms in which each agent updates all decision variables.
In contrast, we consider agents with a common objective function and require that
each agent update only a small subset of decision variables. This has
the advantage that an individual agent's computational burden can be small, even
when solving high-dimensional problems.

The rest of the paper is organized as follows. Section~\ref{sec:prob} includes our problem statement, assumptions, and algorithm. Section~\ref{sec:prelim} provides background
on hybrid systems. We present our hybrid system model in Section~\ref{sec:hsm} and establish the existence of complete solutions. Section~\ref{sec:convergence} proves
that the hybrid multi-agent update law exponentially converges to the minimizer of an objective function. 
We include numerical results as validation in Section~\ref{sec:sim}. 

\section{Problem Statement and Algorithm Overview} \label{sec:prob}
In this section, we state the class of problems that we consider and give an overview of the proposed hybrid optimization algorithm. 

\subsection{Problem Formulation}
We consider a group of~$N$ agents jointly solving an optimization problem of the following form:
\begin{problem} \label{prob:first}
Given an objective function~$L : \R^n \to \R$,
\begin{equation}
\textnormal{minimize } L(x), \qquad x \in \R^n
\end{equation}
using~$N$ distributed agents while requiring that (i) only one agent updates any entry of the decision variable~$x$, and (ii) agents
require only intermittent information sharing from others. 
\end{problem}

Each agent executes computations locally and then shares the results of those computations. Criterion~(i) is there for scalability, only a single agent will update each decision variable. This reduces the computation load on agents and removes duplicated efforts. Criterion~(ii) ensures that the algorithm performs even in environments where communications may be limited. In many practical settings, we expect bandwidth to be limited and/or agents to
have limited onboard power available, which means communications should not be constant.

We assume the following about the objective function~$L$.

\begin{assumption} \label{as:L}
The function~$L$ is twice continuously
differentiable,~$\beta$-strongly convex for some~$\beta > 0$, and~$K$-smooth (namely,~$\nabla L$ is~$K$-Lipschitz). 
\hfill $\triangle$
\end{assumption}

Assumption~\ref{as:L} allows a large number of convex problems to be considered, such as strongly convex quadratic programs. It is a standard assumption in multi-agent 
optimization~\cite{bertsekas89}. 
It implies that~$K \geq \beta$. 

We solve Problem~\ref{prob:first} by applying gradient descent in continuous time using data received intermittently in discrete time. The proposed hybrid optimization algorithm uses jumps to characterize the discrete-time communication events and flows to represent the continuous-time dynamics.
Analogously to past research that has developed distributed versions of the discrete-time gradient
descent law, our update law during flows is based on the following 
(centralized) 
first-order dynamical system:
\begin{equation} \label{eq:gdall}
    \dot{x} + \nabla L(x) =0.
\end{equation}
%Not only are gradient-based controllers already widely used in multi-agent systems, they also provide robustness and are simple to distribute. 
This is motivated by the use of gradient-based controllers in multi-agent systems, e.g., in consensus~\cite{olfati07}, as well
as the simplicity of distributing gradient-based updates and the robustness 
to asynchrony that results from doing so~\cite{bertsekas89}. 
Next, we distribute this across a team of agents.

\subsection{Algorithmic Framework}
We seek to distribute~$\eqref{eq:gdall}$ across a team of agents in
accordance with the parallelization requirement 
in Problem~\ref{prob:first}. 
We consider~$N$ agents indexed over~$i \in [N]:=\lbrace 1, \ldots , N\rbrace$ and divide~$x \in \R^n$ into~$N$ blocks. Then agent~$i$ is responsible for updating and 
communicating values of the~$i$-th block,~$x_{i} \in \mathbb{R}^{n_i}$, 
where~$n_i \in \mathbb{N}$ and~$\sum_{i \in [N]} n_i = n$. 
Thus, the variable~$x$ may be written as the vertical concatenation of all agents' blocks. Each agent performs gradient descent on their own block during flows
but does not update any others. 

Agents' updates occur in continuous time while communication of these updates occurs in discrete time. Communications 
are coordinated using a decreasing timer,~$\tau$, that is shared by all agents. 
When the timer reaches zero, all agents communicate their current values to all of the other agents and the timer resets to a value within a specified interval~$[\tau_{min}, \tau_{max}]$. We assume that communicated data are received at the same time they are sent. 
These communicated blocks are gathered into the vector~$\eta \in \mathbb{R}^n$
with the current value of~$x_i$ being assigned to~$\eta_i$ at communication events.
The value of~$\eta$ is used in each agent's continuous-time computations in a sample-and-hold
manner between communication events. That is, each agent uses the previously communicated data in their updates rather than the continuously evolving values of the other agents. Formally, we write~$\nabla_i L = \frac{\partial L}{\partial x_i}$,
and during flows agent~$i$ executes
\begin{equation}
\dot{x}_i = -\nabla_i L(\eta). 
\end{equation}
This sample-and-hold method is common in the literature~\cite{PHILLIPS2019} and
is used to demonstrate the feasibility of the hybrid approach in multi-agent
optimization.

The complete algorithm is summarized in Algorithm~\ref{alg:first}. 

\begin{algorithm} \label{alg:first}
\SetAlgoLined
Initialization: set~$x_o, \eta_o \in \R^n$ and~$\tau_o \in [0, \tau_{max}]$\;
\While{$\tau \geq 0$}{
$\dot{x}_i = - \nabla_i L(\eta)$, for all~$i\in\lbrace 1, \dots, N \rbrace $\;
$\dot{\tau} = -1$\;
\If{$\tau=0$}{
reset~$\eta_i$ to $x_i$, for all~$i\in\lbrace 1, \dots, N \rbrace $\;
reset~$\tau$ to a value in~$ [\tau_{min},\tau_{max}]$\;}
}
\caption{Distributed Hybrid Gradient Descent}
\end{algorithm}

The next section provides the tools that will be
used to analyze Algorithm~\ref{alg:first}. 

\section{Hybrid System Preliminaries} \label{sec:prelim}
In this section, we recount the background material necessary for the hybrid system
modeling and analysis in the remainder of the paper. 

\subsection{Preliminaries on Hybrid Systems}
For the purposes of this paper, a hybrid system~$\mathcal{H}$ has data~$(C,f,D,G)$ that takes the general form
\begin{align}
    \mathcal{H} = \begin{cases}
    \dot{x} = f(x) &x\in C \\
    x^+ \in G(x) &x\in D
    \end{cases}, \label{eq:hybriddef}
\end{align}
where~$x\in \R^n$ is the system's state, and~$f$ defines the flow map and continuous dynamics for which~$C$ is the flow set. The 
set-valued 
jump map~$G$ captures the system's discrete behavior for the jump set~$D$. More information on this definition and hybrid systems can be found in~\cite{goebel12}. 
\begin{definition}[Hybrid Basic Conditions,~\cite{goebel12}] \label{def:hybridbasic} A hybrid system~$\mathcal{H}$ as in~\eqref{eq:hybriddef} with data~$(C,f,D,G)$ satisfies the~\emph{hybrid basic conditions} if
\begin{itemize}
    \item~$C$ and~$D$ are closed subsets of~$\R^n$;
    \item~$f$ is a continuous function from~$\R^n \to \R^n$;
    \item~$G:\R^n \rightrightarrows \R^n$ is outer semicontinuous and locally bounded relative to~$D$, and~$D \subset \textnormal{dom } G$.
\end{itemize}
\end{definition}
If a hybrid system meets the hybrid basic conditions, then we say that the system is~\emph{well-posed} (Theorem 6.30,~\cite{goebel12}).

We denote solutions to~$\mathcal{H}$ by~$\phi$, which we parameterize by~$(t,j) \in \R_{\geq 0} \times \N$, where~$t$ denotes the ordinary (continuous) time, and~$j$ denotes the jump (discrete) time. Per Definition 2.3 in~\cite{goebel12},~$\textnormal{dom } \phi \subset \R_{\geq 0} \times \N$ is a~\emph{hybrid time domain} if for all~$(T,J) \in \textnormal{dom } \phi$, the~$\textnormal{dom } \phi \cap ([0,T]\times \{0,1,\dots,J\}$ can be written as~$\bigcup_{j=0}^{J-1} ([t_j,t_{j+1}],j)$ for some finite sequence of times~$0 = t_0 \leq t_1 \leq \dots \leq t_J$. We say that a solution~$\phi$ is~\emph{complete} if~$\textnormal{dom } \phi$ is unbounded. A solution~$\phi$ to~$\mathcal{H}$ is called maximal if it cannot be extended further. 

\section{Hybrid System Model} \label{sec:hsm}
In this section, we define a hybrid system model that encompasses all agents' current states and their most recently communicated state values. 
Towards defining this model, we first formally define the timer that governs communication events. This allows us to define the hybrid subsystems that are distributed across the agents. Building on this, we present a definition of the hybrid system modelling the~$N$ agents, their algorithm, and the mechanism governing the communication events. Finally, we show the existence of solutions and conclude that all maximal solutions are complete.

\subsection{Mechanism Governing the Communication Events}
We seek to account for intermittent communication events that occur only at some time instances~$t_j$, for~$j\in \N$, that are not known \emph{a priori}. We assume that the sequence~$\lbrace t_j\rbrace_{j=1}^\infty$ is strictly increasing and unbounded. Between consecutive time events, some amount of time elapses which we upper and lower bound with positive scalars~$\tau_{min}$ and~$\tau_{max}$:
\begin{align}
    0<\tau_{min} \leq t_{j+1} - t_j \leq \tau_{max} \quad \forall j \in \N \setminus \{ 0 \}. \label{eq:timer}
\end{align}
The upper bound~$\tau_{max}$ prevents infinitely long communication delays and ensures convergence, while the lower bound~$\tau_{min}$ rules out Zeno behavior.

To generate events at times~$t_j$ satisfying~\eqref{eq:timer}, let~$\tau$ be the timer that governs when agents exchange data, where~$\tau$ is defined by
\begin{align}
    \dot{\tau} &= -1   &&\tau \in [0, \tau_{max}], \\
    \tau^+ &\in [ \tau_{min}, \tau_{max} ] &&\tau = 0 ,
\end{align}
for~$\tau_{min}, \tau_{max} \in \R_{> 0}$. The timer~$\tau$ steadily decreases until it reaches zero. At this point, it is reset to a value within~$[ \tau_{min}, \tau_{max} ]$.

There is indeterminacy built into the timer in that the reset
map is only confined to a compact interval,~$[\tau_{min}, \tau_{max}]$, where~$\tau_{min}$ and~$\tau_{max}$ are both positive real numbers.

\subsection{Hybrid Subsystems} 
Recall that agent~$i$ stores and updates its own state variable~$x_i \in \R^{n_i}$, and the variable~$x$ is the vertical concatenation of all agents' states. Data from all agents are collectively stored in~$\eta \in \R^{n}$ at communication events. We define the state of agent~$i$'s hybrid system as~$\xi_i = (x_{i}, \eta, \tau)$, where~$x_{i}$ is agent~$i$'s state (the one it is responsible for updating),~$\eta$ is the memory state storing the states of the agents measured at communication events, and~$\tau$ is defined as above. 
This leads to the hybrid subsystem given by
\begin{align}
    \dot{\xi_i} & = \begin{bmatrix}
    -\nabla_{i} L(\eta) \\
    0 \\
    -1
    \end{bmatrix} &&\xi_i \in \R^{n_i} \times \R^n \times [0,\tau_{max}]  \\
    \xi_i^+ &\in \begin{bmatrix}
    x_{i} \\
    x \\
    [\tau_{min}, \tau_{max}]
    \end{bmatrix} &&\xi_i \in \R^{n_i} \times \R^n \times \{0\},
\end{align}
where~$x=(x_1^T, \ldots, x_N^T)^T$.

\subsection{Combined Hybrid System}
We are now ready to combine the distributed subsystems into one hybrid system for analysis. First, we define a variable~$z = (z_1, z_2) \in \R^n \times \R^n$ such that
\begin{align}
    z_1 &= \col (x_1, \dots, x_N) \\
    z_2 &= \eta,
\end{align}
where~$\col (x_1, \dots, x_N) = (x_1^T, \ldots, x_N^T)^T$. 

We define the state of the combined hybrid system as~$\xi = (z_1, z_2, \tau) \in \mathcal{X}$, where~$z_1$,~$z_2$, and~$\tau$ are defined as above, and~$\mathcal{X}:=\R^n \times \R^n \times [0,\tau_{max}]$.
This leads to the combined hybrid system~$\mathcal{H}=(C,f,D,G)$ given by
\begin{align}
    \dot{\xi} & = \begin{bmatrix}
    -\nabla L(z_2) \\
    0 \\
    -1
    \end{bmatrix}:= f(\xi)  && \xi \in C, \label{eq:hfirst} \\ 
    \xi^+ &\in \begin{bmatrix}
    z_1 \\
    z_1 \\
    [\tau_{min}, \tau_{max}]
    \end{bmatrix}:= G(\xi) && \xi \in D, \label{eq:hlast} 
\end{align}
where~$C:= \mathcal{X}$ and~$D:=\R^n \times \R^n \times \lbrace 0 \rbrace$.
 
\subsection{Hybrid Basic Conditions }
We now demonstrate that~$\mathcal{H}$ meets the hybrid basic conditions and is well-posed.
\begin{lemma}
Let~$L$ satisfy Assumption~\ref{as:L}. Then, the hybrid system given by~$\mathcal{H}$ with data~$(C,f,D,G)$ defined in~\eqref{eq:hfirst}--\eqref{eq:hlast} satisfies the hybrid basic conditions from Definition~\ref{def:hybridbasic} and is nominally well-posed as a result.
\end{lemma}
\noindent \emph{Proof: } \ifbool{Report}{The sets~$C$ and~$D$ are closed subsets of~$\R^n \times \R^n \times \R$ by definition. Due to our assumption that~$\nabla L$ is continuous,~$f$ is a continuous function from~$\R^n \times \R^n \times \R$ to~$\R^n \times \R^n \times \R$. By construction,~$G$ is outer semicontinuous and locally bounded relative to~$D$. Finally,~$D \subset \textnormal{dom } G$ because~$\textnormal{dom } G$ is~$\R^n \times \R^n \times \R$. \hfill $\blacksquare$}{See~\cite{hendrickson21}.\hfill $\blacksquare$}

%The sets~$C$ and~$D$ are closed subsets of~$\R^n \times \R^n \times \R$ by definition. Due to our assumption that~$\nabla L$ is continuous,~$f$ is a continuous function from~$\R^n \times \R^n \times \R$ to~$\R^n \times \R^n \times \R$. By construction,~$G$ is outer semicontinuous and locally bounded relative to~$D$. Finally,~$D \subset \textnormal{dom } G$ because~$\textnormal{dom } G$ is~$\R^n \times \R^n \times \R$. \hfill $\blacksquare$

\subsection{Existence of Solutions}
In addition to being well-posed, there exists a nontrivial solution to~$\mathcal{H}$ from each point in~$C\cup D$, and all maximal solutions are complete and not Zeno under mild conditions on problem
parameters. Complete solutions cannot be extended further and their domains are unbounded. Practically, this means that the proposed algorithm may run for an arbitrarily long period of time and does not reach a point where it can neither flow nor jump.
\begin{lemma}[Existence of Solutions] \label{lem:solutions}
Let Assumption~\ref{as:L} hold. Let~$\tau_{min}$ and~$\tau_{max}$ be such that~$0 < \tau_{min} \leq \tau_{max} < \frac{\beta^2}{3K^3}$, where~$\beta$ is the strong convexity constant of~$L$ and~$K$ is the Lipschitz constant of~$\nabla L$. 
Then there exists a nontrivial solution to~$\mathcal{H}=(C,f,D,G)$ from every initial point in~$C \cup D$. Additionally, every maximal solution~$\phi$ to the hybrid system~$\mathcal{H}$ is complete and not Zeno.
\end{lemma}
\noindent \emph{Proof: } \ifbool{Report}{See the appendix. \hfill $\blacksquare$}{See~\cite{hendrickson21}.\hfill $\blacksquare$}

\section{Convergence Analysis} \label{sec:convergence}
In this section, we define the set for solutions to converge to and present some useful properties of the hybrid system~$\mathcal{H}$ in Lemmas~\ref{lem:contract} and~\ref{lem:gradprop}. We then propose a Lyapunov function in Lemma~\ref{lem:comparison}. As an interim result, we show that for a solution~$\phi = (\phi_{z_1},\phi_{z_2},\phi_{\tau})$ to~$\mathcal{H}$ in~\eqref{eq:hfirst}--\eqref{eq:hlast}, if~$\phi_{z_1}(0,0) = \phi_{z_2}(0,0)$, we are able to bound the distance from the minimizer of~$L$ for all~$(t,j) \in \textnormal{dom }\phi$. Finally, we present our main result, exponential convergence to the minimizer of~$L$, in Theorem~\ref{thm:ges}.

\subsection{Convergence Set}
 Let~$\textbf{0}_n$ be the vector of zeros in~$\R^n$; similarly, let~$\textbf{0}_{n_i}$ be the vector of zeros in~$\R^{n_i}$. Convergence using gradient descent occurs when the gradient of~$L$ is~$\textbf{0}_n$. Given a complete solution~$\phi = (\phi_{z_1},\phi_{z_2},\phi_{\tau})$ to the hybrid system~$\mathcal{H}$, we seek to assure that~$\lim_{t+j \to \infty} \nabla_i L(\phi_{z_2}(t,j)) = \textbf{0}_{n_i},$ for~$i=1, \dots, N$. This is equivalent to a set convergence problem where the set to converge to for the hybrid system~$\mathcal{H}$ is given by \ifbool{Report}{
 \begin{align}
    \mathcal{A} &:= \lbrace \xi = (z_1, z_2, \tau) \in \mathcal{X} : \nabla L(z_2) = \textbf{0}_n, z_2 = z_1, \tau \in [0, \tau_{max}] \rbrace \\
    & =\{x^*\}\times\{x^*\}\times[0, \tau_{max}], \label{eq:Adef}
\end{align}}{
\begin{align}
    \mathcal{A} &:= \lbrace \xi = (z_1, z_2, \tau) \in \mathcal{X} : \\
    &\quad \quad \nabla L(z_2) = \textbf{0}_n, z_2 = z_1, \tau \in [0, \tau_{max}] \rbrace \\
    & =\{x^*\}\times\{x^*\}\times[0, \tau_{max}], \label{eq:Adef}
\end{align}}
where~$x^*$ is the unique fixed point of~$\nabla L$. Equivalence of the expression for~$\mathcal{A}$ stems from Assumption~\ref{as:L}: because~$L$ is strongly convex, it has a unique minimum (denoted by~$x^*$) and this unique minimum is the unique stationary point of~$\nabla L$. Given a vector~$\xi=(z_1,z_2, \tau) \in \mathcal{X}$, the squared distance from~$\mathcal{A}$ is given by~$|\xi|^2_\mathcal{A} := \|z-z^*\|^2 = \|z_1 - x^*\|^2 + \|z_2 - x^*\|^2$, where~$\|\cdot\|$ denotes the Euclidean norm throughout this paper.

\subsection{Useful Properties of~$\mathcal{H}$}
Combining gradient descent  with a bound on~$\tau_{\max}$ allows us to establish relationships that prove useful during Lyapunov analysis.

\begin{lemma} \label{lem:contract}
Let Assumption~\ref{as:L} hold. Consider the hybrid system given by~$\mathcal{H}$ with data~$(C,f,D,G)$ defined
in~\eqref{eq:hfirst}--\eqref{eq:hlast}. 
Let~$\tau_{min}$ and~$\tau_{max}$ be such that~$0 < \tau_{min} \leq \tau_{max} < \frac{\beta^2}{3K^3}$, where~$\beta$ is the strong convexity constant of~$L$ and~$K$ is the Lipschitz constant of~$\nabla L$. Denote the unique fixed point of~$\nabla L$ by~$x^*$. Pick a solution~$\phi = (\phi_{z_1},\phi_{z_2},\phi_{\tau})$ to~$\mathcal{H}$ such that $\phi_{z_1}(0,0) = \phi_{z_2}(0,0)$. For each~$I^j := \{ t: (t,j) \in \textnormal{dom } \phi \}$ with nonempty interior and with~$t_{j+1} > t_j$ such that~$[t_j, t_{j+1}] = I^j$, we have
 \begin{align}
     \phi_{z_1}(t,j) &= \phi_{z_2}(t_j,j) - (t-t_j) \nabla L(\phi_{z_2}(t_j,j)) \label{eq:z1update} \\
     \phi_{z_2}(t,j) &= \phi_{z_2}(t_j,j), \label{eq:z2update}
 \end{align}
for all~$t \in (t_j, t_{j+1})$.
 Additionally, for all~$(t,j) \in \textnormal{dom } \phi$, the following are satisfied:
 \begin{align}
     \|\phi_{z_1}(t,j) - x^* \|^2 &\leq q(t,t_j) \|\phi_{z_2}(t_j,j) - x^*\|^2; \quad \label{lem:conpart1}\\
     \|\phi_{z_1}(t,j) - \phi_{z_2}(t,j) \| &\leq \tau_{max} \| \nabla L(\phi_{z_2}(t_j,j)) \|; \label{lem:conpart2} \\
     \|\phi_{z_1}(t,j) - x^* \|^2 &\geq B \|\phi_{z_2}(t_j,j) - x^*\|^2; \label{lem:conpart3}
 \end{align}
 where~$q (t,t_j):= (1 - 2(t-t_j) \beta + (t-t_j)^2 K^2) \in (0,1)$ and~$B:=(1-2\tau_{max} K) \in (0,1)$.
\end{lemma}
\noindent \emph{Proof: } \ifbool{Report}{ Given~$t\in(t_j, t_{j+1})$, the solution~$\phi$ has flowed some distance given by~$(t-t_j)\dot{\phi}$ 
where~$\dot{\phi}$ is constant due to the sample-and-hold
methodology.
Applying our definition of~$f$ in~\eqref{eq:hfirst} gives~\eqref{eq:z1update} 
and~\eqref{eq:z2update}. \\
\emph{Proof of}~\eqref{lem:conpart1}: Using~\eqref{eq:z1update} and the fact that~$\nabla L(x^*) = 0$, we can rewrite~$\|\phi_{z_1}(t,j) - x^*\|^2$ as
\begin{align}
    \|\phi_{z_1}(t,j) - x^*\|^2 &= \|\phi_{z_2}(t_j,j) - (t-t_j) \nabla L(\phi_{z_2}(t_j,j)) - x^* + (t-t_j) \nabla L(x^*)\|^2 \\
    &= \|\phi_{z_2}(t_j,j) - x^*\|^2 -2(t\!-\!t_j) (\nabla L(\phi_{z_2}(t_j,j))\!-\!\nabla L(x^*))^T (\phi_{z_2}(t_j,j)\!-\!x^*) \\
    &\quad + (t-t_j)^2 \|\nabla L(\phi_{z_2}(t_j,j)) - \nabla L(x^*)\|^2,
\end{align}
where the second equality follows from expanding the norm squared. Using the~$\beta$-strong convexity of~$L$ and the Lipschitz property of~$\nabla L$, we upper bound this with
\begin{align}
    \|\phi_{z_1}(t,j) - x^*\|^2 &\leq \|\phi_{z_2}(t_j,j) - x^*\|^2\! -\!2 (t\!-\!t_j) \beta \|\phi_{z_2}(t_j,j) - x^*\|^2 + (t-t_j)^2 K^2 \|\phi_{z_2}(t_j,j) - x^*\|^2 \\
    &= (1\!-\!2(t\!-\!t_j) \beta\!+\!(t\!-\!t_j)^2 K^2) \|\phi_{z_2}(t_j,j)\!-\!x^*\|^2.
\end{align}

For contraction, we must show~$(1\!-\!2(t\!-\!t_j) \beta\!+\!(t\!-\!t_j)^2 K^2) < 1$.
To derive a sufficient condition for this, 
note that~$(1\!-\!2(t\!-\!t_j) \beta\!+\!(t\!-\!t_j)^2 K^2) < 1$ 
may be rewritten as~$(t-t_j) K^2 < 2 \beta$ by subtracting~$1$ from both sides, dividing by~$t-t_j$, and then adding~$2\beta$ to both sides. Using~$\beta \leq K$ and~$\tau_{max} < \frac{\beta^2}{3K^3}$, we have~$(t-t_j) \leq \tau_{max} < \frac{\beta^2}{3K^3} < \frac{2\beta}{K^2}$, and therefore~$(1\!-\!2(t\!-\!t_j) \beta\!+\!(t\!-\!t_j)^2 K^2) < 1$. To show that this term is also positive, it is sufficient to show that~$(1 - 2(t-t_j) \beta) \geq 1 - 2 \tau_{max} K > 0$. This is satisfied for $\tau_{max} < \frac{\beta^2 }{3 K^3} \leq \frac{1}{3 K} < \frac{1}{K}$. Thus,~$q(t,t_j) = (1 - 2(t-t_j) \beta + (t-t_j)^2 K^2) \in (0,1)$.

\emph{Proof of}~\eqref{lem:conpart2}: From~\eqref{eq:z1update}--\eqref{eq:z2update}, we have~$\| \phi_{z_1}(t,j)\!-\!\phi_{z_2}(t,j) \| = \|\phi_{z_2}(t_j,j)\!-\!(t\!-\!t_j) \nabla L(\phi_{z_2}(t_j,j))\!-\!\phi_{z_2}(t_j,j) \|$. Simplifying and bounding~$t-t_j$ above by~$\tau_{max}$ gives the final result.

\emph{Proof of}~\eqref{lem:conpart3}: Using~\eqref{eq:z1update} and the fact that~$\nabla L(x^*) = 0$, we can rewrite~$\|\phi_{z_1}(t,j) - x^*\|^2$ as $\|\phi_{z_2}(t_j,j)\!-\!(t\!-\!t_j) \nabla L(\phi_{z_2}(t_j,j))\!-\!x^*\!+\!(t\!-\!t_j) \nabla L(x^*)\|^2$. Expanding the norm squared gives
\begin{align}
    \|\phi_{z_1}(t,j) - x^*\|^2 &= \|\phi_{z_2}(t_j,j)\!-\!x^*\|^2\! +(t\!-\!t_j)^2 \|\nabla L(\phi_{z_2}(t_j,j))\!-\!\nabla L(x^*)\|^2\\
    &\quad \!-\!2(t\!-\!t_j)(\phi_{z_2}(t_j,j)\!-\!x^*)^T(\nabla L(\phi_{z_2}(t_j,j))\!-\!\nabla L(x^*)). 
\end{align}
Dropping the middle term (which is positive) 
and using~$(\phi_{z_2}(t_j,j) - x^*)^T(\nabla L(\phi_{z_2}(t_j,j))- \nabla L(x^*)) \leq K \|\phi_{z_2}(t_j,j) - x^*\|^2$, which follows from~$\nabla L$ being~$K$-Lipschitz, we can derive a lower bound:
\begin{align}
    \|\phi_{z_1}(t,j) - x^* \|^2 &\geq (1 -2(t-t_j) K) \|\phi_{z_2}(t_j,j) - x^*\|^2\\
    &\geq (1 -2\tau_{max} K) \|\phi_{z_2}(t_j,j) - x^*\|^2. \label{eq:z1geqz2}
\end{align}
Let~$B := 1 - 2\tau_{max} K$. Then using~$\beta \leq K$, we have $\tau_{max} < \frac{\beta^2 }{6 K^3}  < \frac{1}{2K}$. Therefore,~$B \in (0,1)$.
\hfill $\blacksquare$}{See~\cite{hendrickson21}. \hfill $\blacksquare$}

In preparation for establishing the convergence
properties of~$\mathcal{H}$, we also show that the angle between the gradient of the current state and the gradient of the previously communicated state is never greater than~$90$ degrees as a result of the bound on~$\tau_{max}$. This is formally stated in Lemma~\ref{lem:gradprop}. 
\begin{lemma} \label{lem:gradprop}
Let Assumption~\ref{as:L} hold. Consider the hybrid system given by~$\mathcal{H}$ with data~$(C,f,D,G)$ defined
in~\eqref{eq:hfirst}--\eqref{eq:hlast}. Let~$\tau_{min}$ and~$\tau_{max}$ be such that~$0 < \tau_{min} \leq \tau_{max} < \frac{\beta^2}{3K^3}$, where~$\beta$ is the strong convexity constant of~$L$ and~$K$ is the Lipschitz constant of~$\nabla L$. Denote the unique fixed point of~$\nabla L$ by~$x^*$. Pick a solution~$\phi$ such that~$\phi_{z_1}(0,0) = \phi_{z_2}(0,0)$. For each~$I^j := \{ t: (t,j) \in \textnormal{dom } \phi \}$ with nonempty interior and with~$t_{j+1} > t_j$ such that~$[t_j, t_{j+1}] = I^j$, we have
\begin{align}
    \nabla L(\phi_{z_1}(t,j))^T \nabla L(\phi_{z_2}(t,j)) \geq A \|\phi_{z_2}(t_j,j)-x^*\|^2,
\end{align}
for all~$t \in (t_j, t_{j+1})$, where~$A:=\beta^2(1 - 2 \tau_{max} K) - \tau_{max} K^3 > 0$.
\end{lemma}
\noindent \emph{Proof: }\ifbool{Report}{We begin by expanding~$\|\nabla L(\phi_{z_1}(t,j))\|^2$:
\begin{align}
    \|\nabla L(\phi_{z_1}(t,j))\|^2 &= \nabla L(\phi_{z_1}(t,j))^T\Big(\nabla L(\phi_{z_2}(t,j)) + \nabla L(\phi_{z_1}(t,j)) - \nabla L(\phi_{z_2}(t,j))\Big) \\
    &\leq \nabla L(\phi_{z_1}(t,j))^T \nabla L(\phi_{z_2}(t,j)) + \|\nabla L(\phi_{z_1}(t,j))\| \|\nabla L(\phi_{z_1}(t,j)) - \nabla L(\phi_{z_2}(t,j))\| \\
    &\leq \nabla L(\phi_{z_1}(t,j))^T \nabla L(\phi_{z_2}(t,j)) + K \|\nabla L(\phi_{z_1}(t,j))\| \|\phi_{z_1}(t,j) - \phi_{z_2}(t,j)\| \\
    &\leq \nabla L(\phi_{z_1}(t,j))^T \nabla L(\phi_{z_2}(t,j)) + \tau_{max} K \|\nabla L(\phi_{z_1}(t,j))\| \|\nabla L(\phi_{z_2}(t_j,j))\|, \label{eq:z1exp}
\end{align}
where the last two inequalities are from the~$K$-Lipschitz property of~$\nabla L$ which states $\|\nabla L(x) - \nabla L(y)\| \leq K \|x-y\|,$ for all~$x,y \in \R^n$ and~\eqref{lem:conpart2} in Lemma~\ref{lem:contract}. Note that the~$K$-Lipschitz property also gives the inequalities~$\|\nabla L(\phi_{z_1}(t,j))\| \leq K \|\phi_{z_1}(t,j)-x^*\|$ and~$\|\nabla L(\phi_{z_2}(t_j,j))\| \leq K \|\phi_{z_2}(t_j,j)-x^*\|$ because~$\nabla L(x^*) = 0$. Thus,~\eqref{eq:z1exp} becomes
\begin{align}
    \|\nabla L(\phi_{z_1}(t,j))\|^2 &\leq \nabla L(\phi_{z_1}(t,j))^T \nabla L(\phi_{z_2}(t_j,j)) + \tau_{max} K^3 \|\phi_{z_1}(t,j)-x^*\| \|\phi_{z_2}(t_j,j)-x^*\| \\
    &\leq \nabla L(\phi_{z_1}(t,j))^T \nabla L(\phi_{z_2}(t_j,j)) + \tau_{max} K^3 \|\phi_{z_2}(t_j,j)-x^*\|^2, \label{eq:z1Upper}
\end{align}
where the last inequality follows from~\eqref{lem:conpart1} in Lemma~\ref{lem:contract}. Because~$L$ is~$\beta$-strongly convex, we have~$\|\nabla L(\phi_{z_1}(t,j)) - \nabla L(x^*) \|^2 \geq \beta^2 \|\phi_{z_1}(t,j) - z^*\|^2$. To lower bound~\eqref{eq:z1Upper} in terms of~$\|\phi_{z_2}(t_j,j) - x^*\|^2$, we combine this with~$\nabla L(x^*) = 0$ and~\eqref{lem:conpart3} in Lemma~\ref{lem:contract}. This gives the inequality
\begin{align}
\|\nabla L(\phi_{z_1}(t,j))\|^2 &\geq \beta^2 \|\phi_{z_1}(t,j)-x^*\|^2 \\
&\geq \beta^2 B \|\phi_{z_2}(t_j,j) - x^*\|^2, \label{eq:z1Lower2}
\end{align}
where~$B:=(1 - 2 \tau_{max} K)$. Then using~\eqref{eq:z1Lower2} in conjunction with~\eqref{eq:z1Upper} allows us to lower bound~$\nabla L(\phi_{z_1}(t,j))^T \nabla L(\phi_{z_2}(t,j))$:
\begin{align}
    \nabla L(\phi_{z_1}(t,j))^T \nabla L(\phi_{z_2}(t,j)) &\geq \beta^2 B \|\phi_{z_2}(t_j,j) - x^*\|^2  - \tau_{max} K^3 \|\phi_{z_2}(t_j,j)-x^*\|^2 \\
    &=(\beta^2(1 - 2 \tau_{max} K) - \tau_{max} K^3) \|\phi_{z_2}(t_j,j)-x^*\|^2 .
\end{align}
This lower bound is positive for~$\tau_{max} < \frac{\beta^2 }{2 \beta^2K + K^3}$.
Because~$\beta \leq K$, the denominator~$2\beta^2K + K^3 \leq 3K^3$. This lower bound is positive since~$\tau_{max} < \frac{\beta^2 }{3 K^3} \leq \frac{\beta^2 }{2 \beta^2K + K^3}$. \hfill $\blacksquare$}{See~\cite{hendrickson21}. \hfill $\blacksquare$}

\subsection{Bound on the Lyapunov Function}
Central to proving our main result is a Lyapunov function that is bounded above and below by~$\mathcal{K}_\infty$ comparison functions~$\alpha_1, \alpha_2$ given in Lemma~\ref{lem:comparison}. 
%These bounds hold when Assumption~\ref{as:L} holds and~$\tau_{max}$ is bounded but do not impose any requirements on initialization.

\begin{lemma} \label{lem:comparison}
Let Assumption~\ref{as:L} hold. Let~$\tau_{min}$ and~$\tau_{max}$ be such that~$0 < \tau_{min} \leq \tau_{max} < \frac{\beta^2}{3K^3}$, where~$\beta$ is the strong convexity constant of~$L$ and~$K$ is the Lipschitz constant of~$\nabla L$. Let~$V : \mathcal{X} \to \R_{\geq 0}$ be a Lyapunov function candidate for the hybrid system~$\mathcal{H}= (C,f,D,G)$ defined 
in~\eqref{eq:hfirst}--\eqref{eq:hlast}, given by
\begin{align}
    V(\xi) &= (L(z_1) - L(x^*))^2 + (L(z_2) - L(x^*))^2,
\end{align}
for all~$\xi=(z_1, z_2, \tau) \in \mathcal{X}$, where~$L$ is the objective function and~$x^*$ is the unique fixed point of~$\nabla L$. Then there exist~$\alpha_1, \alpha_2 \in \mathcal{K}_\infty$ such that
\begin{align}
    &\alpha_1(|\xi|_\mathcal{A}) \leq V(\xi) \leq \alpha_2(|\xi|_\mathcal{A})
\end{align}
for all~$\xi \in C \cup D \cup G(D)$. In particular,~$\alpha_1$ and~$\alpha_2$ may be given by, for each~$s\geq 0$,  
\begin{align}
    \alpha_1(s) = \frac{\beta^2}{16} s^4 \quad \textnormal{and} \quad \alpha_2(s) = \frac{K^2}{2} s^4.
\end{align}
\end{lemma}
\noindent \emph{Proof: } \ifbool{Report}{The minimizer of~$L$ is~$\xi^* = (z_1^*, z_2^*, \tau) = (x^*, x^*, \tau)$ 
for any~$\tau$, 
and, by construction,~$V(\xi)$ is zero only for~$\xi = \xi^*$ and is positive otherwise.

Because~$\nabla L$ is~$K$-Lipschitz,~$L(x)-L(x^*) \leq \nabla L(x^*)^T (x-x^*) + \frac{K}{2} \|x-x^*\|^2 =  \frac{K}{2} \|x-x^*\|^2,$ for all~$x\in \R^n$. Thus,~$V(\xi)$ may be bounded as
\begin{align}
    V(\xi) \leq \frac{K^2}{4} \|z_1 - x^*\|^4 + \frac{K^2}{4} \|z_2 - x^*\|^4. \label{eq:lxi}
\end{align}
First consider the case where~$\|z_1 - x^*\| \leq \|z_2 - x^*\|$. Then~\eqref{eq:lxi} becomes
\begin{align}
    V(\xi) \leq \frac{K^2}{2} \|z_2 - x^*\|^4 \leq \frac{K^2}{2} |\xi|^4_\mathcal{A}, \label{eq:compupper}
\end{align}
where the last inequality follows from~$\|z_2 - x^*\|^2 \leq |\xi|_\mathcal{A}^2$ by definition of~$|\xi|_\mathcal{A}$. If~$\|z_1 - x^*\| > \|z_2 - x^*\|$ instead, then~\eqref{eq:lxi} is bounded above by
\begin{align}
    V(\xi) < \frac{K^2}{2} \|z_1 - x^*\|^4 \leq \frac{K^2}{2} |\xi|^4_\mathcal{A},
\end{align}
where the definition of~$|\xi|_\mathcal{A}$ is again used in the last inequality. Thus,~$V(\xi) \leq \frac{K^2}{2} |\xi|^4_\mathcal{A}$ for all~$\xi \in \mathcal{X}$ and we set~$\alpha_2(s) = \frac{K^2}{2} s^4 \in \mathcal{K}_\infty$ for all~$s \geq 0$.

The~$\beta$-strong convexity of~$L$ and~$\nabla L(x^*) = 0$ allow us to write~$L(z_1) - L(x^*) \geq \frac{\beta}{2} \|z_1-x^*\|^2$ 
and~$L(z_2) - L(x^*) \geq \frac{\beta}{2} \|z_2-x^*\|^2$.
Applying both inequalities to the definition of~$V$ gives
\begin{align}
    V(\xi) &\geq \frac{\beta^2}{4} \|z_1-x^*\|^4 + \frac{\beta^2}{4} \|z_2-x^*\|^4. \label{eq:vupper}
\end{align}
First, consider the case where~$\|z_1 - x^*\| \geq \|z_2 - x^*\|$. Then dropping the second term in~\eqref{eq:vupper} and using~$|\xi|^2_\mathcal{A} = \|z_1 - x^*\|^2 + \|z_2 - x^*\|^2 \leq 2 \|z_1 - x^*\|^2$,  
\begin{align}
    V(\xi) &\geq \frac{\beta^2}{4}\|z_1-x^*\|^4 \geq \frac{\beta^2}{16} |\xi|_\mathcal{A}^4. \label{eq:complower}
\end{align}
Now consider the case that~$\|z_1 - x^*\| < \|z_2 - x^*\|$. Then the same steps apply and~$V(\xi) > \frac{\beta^2}{16} |\xi|_\mathcal{A}^4$. Thus,~$V(\xi) \geq \frac{\beta^2}{16} |\xi|_\mathcal{A}^4$ for all~$\xi \in \mathcal{X}$.
Accordingly,~$\alpha_1(s) = \frac{\beta^2}{16} s^4$ for all~$s\geq 0$. \hfill $\blacksquare$}{See~\cite{hendrickson21}.\hfill $\blacksquare$}

\subsection{Exponential Convergence}
Using Lemmas~\ref{lem:contract},~\ref{lem:gradprop}, and~\ref{lem:comparison}, we are able to bound the distance to the minimizer of~$L$ over time for a class of initial conditions in Proposition~\ref{prop:conv}. This result will then be expanded to include all possible solutions and initial conditions in Theorem~\ref{thm:ges}, thus showing exponential convergence to the minimizer of~$L$.   
\begin{proposition} \label{prop:conv}
Let Assumption~\ref{as:L} hold and consider the hybrid system~$\mathcal{H}$ defined in~\eqref{eq:hfirst}-\eqref{eq:hlast}. Let~$\mathcal{A}$ be as defined in~\eqref{eq:Adef} and let~$\tau_{min}$ and~$\tau_{max}$ be such that~$0 < \tau_{min} \leq \tau_{max} < \frac{\beta^2}{3K^3}$, where~$\beta$ is the strong convexity constant of~$L$ and~$K$ is the Lipschitz constant of~$\nabla L$. For each solution~$\phi$ to~$\mathcal{H}$ such that~$\phi_{z_1}(0,0) = \phi_{z_2}(0,0)$, for all~$(t,j) \in \textnormal{dom } \phi$, the following is satisfied:
\begin{align}
    |\phi(t,j)|_\mathcal{A} &\leq \sqrt{\frac{K}{\beta}} \sqrt[4]{8} \exp \Big( -\frac{\beta A B}{8 K^2} t\Big) |\phi(0,0)|_\mathcal{A},
\end{align}
where~$A=\beta^2(1 - 2 \tau_{max} K) - \tau_{max} K^3 > 0$ and~$B=(1-2\tau_{max} K) \in (0,1)$. 
\end{proposition}
\noindent \emph{Proof: } \ifbool{Report}{We first consider~$\xi \in C$ and the Lyapunov function~$V$ defined in Lemma~\ref{lem:comparison}. The partial derivatives of~$V$ with respect to~$z_i$ are given by~$\nabla_{z_i} V(\xi) = 2 \nabla L(z_i) (L(z_i) - L(x^*)) \in \R^n$, for~$i=1,2$, where~$x^*$ is the unique fixed point of~$\nabla L$.
This leads to 
\begin{align}
    \langle \nabla V(\xi), f(\xi) \rangle =\!-\!2(L(z_1)\!-\!L(x^*)) \nabla L(z_1)^T \nabla L(z_2). \quad \label{eq:innerdef}
\end{align}
We now pick a solution~$\phi$ such that~$\phi_{z_1}(0,0) = \phi_{z_2}(0,0)$. For each~$I^j := \{t:(t,j) \in \textnormal{dom } \phi \}$ with nonempty interior and with~$t_{j+1} > t_j$ such that~$[t_j, t_{j+1}] = I^j$, we have from Lemma~\ref{lem:gradprop},
\begin{align}
    \nabla L(\phi_{z_1}(t,j))^T \nabla L(\phi_{z_2}(t,j)) &\geq A \|\phi_{z_2}(t_j,j)-x^*\|^2,
\end{align}
where~$A=\beta^2(1 - 2 \tau_{max} K) - \tau_{max} K^3 >0$. Combining this with~$(L(\phi_{z_1}(t,j)) - L(x^*)) \geq \frac{\beta}{2} \|\phi_{z_1}(t,j) - x^*\|^2$ from the~$\beta$-strong convexity of~$L$, and using~\eqref{eq:innerdef} we can write
\begin{align}
    \langle \nabla V(\phi(t,j)), f(\phi(t,j)) \rangle
    &\leq - \beta\|\phi_{z_1}(t,j) - x^*\|^2 \nabla L(\phi_{z_1}(t,j))^T \nabla L(\phi_{z_2}(t,j)) \\
    &\leq - \beta A \|\phi_{z_1}(t,j) - x^*\|^2  \|\phi_{z_2}(t_j,j) - x^*\|^2 \\
    &\leq - \beta A B \|\phi_{z_2}(t_j,j)-x^*\|^4, \label{eq:innerz2}
\end{align}
where the last inequality is from applying~\eqref{lem:conpart3} in Lemma~\ref{lem:contract}.
Additionally, for each~$(t,j)$ in the interval of flow~$[t_j, t_{j+1}] \times \{j\}$ in~$\textnormal{dom } \phi$, we have~$\|\phi_{z_1}(t,j) - x^*\| \leq \|\phi_{z_2}(t_j,j) - x^*\|$ from~\eqref{lem:conpart1} in Lemma~\ref{lem:contract}. Applying this to the definition of~$|\phi(t,j)|_\mathcal{A}^2$ gives the relationship~$|\phi(t,j)|_\mathcal{A}^2 \leq 2 \|\phi_{z_2}(t_j,j) -x^*\|^2$. Squaring both sides allows us to rewrite~\eqref{eq:innerz2} as
\begin{align}
    \langle \nabla V(\phi(t,j)), f(\phi(t,j)) \rangle &\leq - \frac{\beta A B}{4} |\phi(t,j)|_\mathcal{A}^4. \label{eq:dotVupper}
\end{align}
This bound may be related back to~$V(\phi (t,j))$ using the comparison function~$\alpha_2$ from Lemma~\ref{lem:comparison}, which leads to
\begin{align}
    \langle \nabla V(\phi(t,j)), f(\phi(t,j)) \rangle &\leq - \frac{\beta A B}{2K^2} V(\phi(t,j)). \label{eq:tinbound}
\end{align}

We now consider the change of~$V$ at jumps. For each~$\xi \in D,$ and~$g=(g_{z_1},g_{z_2},g_{\tau}) \in G(\xi)$,~$V(g) = (L(g_{z_1}) - L(x^*))^2 + (L(g_{z_2}) - L(x^*))^2$. Thus, we can write the change of~$V$ at jumps as~$V(g) - V(\xi) = (L(g_{z_1}) - L(x^*))^2 - (L(z_2) - L(x^*))^2$ because~$g_{z_1}=z_1$ and~$g_{z_2} = z_1$ from~\eqref{eq:hlast}.

For a solution~$\phi$ such that~$\phi_{z_1}(0,0) = \phi_{z_2}(0,0)$, this is equivalent to
\begin{align}
    V(G(\phi(t_{j+1},j))) - V(\phi(t_{j+1},j)) = (L(\phi_{z_1}(t_{j\!+\!1},j\!+\!1))\!-\!L(x^*))^2\! -\!(L(\phi_{z_2}(t_{j\!+\!1},j))\!-\!L(x^*))^2,
\end{align}
for all~$(t_{j+1},j), (t_{j+1},j+1) \in \textnormal{dom }\phi$. For this quantity to be nonpositive, it is sufficient to show that~$L(\phi_{z_1}(t_{j+1},j+1)) \leq L(\phi_{z_2}(t_{j+1},j))$.

Towards doing this, we leverage~$L(x) \leq L(y) + \nabla L(y)^T(x - y) + \frac{K}{2}\|x - y\|^2$, for all~$x,y\in \R^n$ from the~$K$-Lipschitz property of~$\nabla L$. Applying this for~$x=\phi_{z_1}(t_{j+1},j+1), y=\phi_{z_2}(t_{j+1},j)$, we have
\begin{align}
L(\phi_{z_1}(t_{j+1},j+1)&) \leq L(\phi_{z_2}(t_{j+1},j)) +\nabla L(\phi_{z_2}(t_{j+1},j))^T(\phi_{z_1}(t_{j+1},j+1)\!-\!\phi_{z_2}(t_{j+1},j)) \\
& \qquad \qquad + \frac{K}{2}\|\phi_{z_1}(t_{j+1},j+1)-\phi_{z_2}(t_{j+1},j)\|^2 \\
&= L(\phi_{z_2}(t_{j+1},j))\!+\!\nabla L(\phi_{z_2}(t_{j+1},j))^T(\phi_{z_2}(t_j,j)\!-\!(t_{j+1}\!-\!t_j) \nabla L(\phi_{z_2}(t_j,j))\!-\!\phi_{z_2}(t_{j+1},j)) \\
& \qquad \qquad   +\frac{K}{2}\|\phi_{z_2}(t_j,j)\!-\!(t_{j+1}\!-\!t_j) \nabla L(\phi_{z_2}(t_j,j))\!-\!\phi_{z_2}(t_j,j)\|^2 \\
&= L(\phi_{z_2}(t_{j+1},j)) - \Big(1- \frac{K}{2}(t_{j+1} - t_j)\Big) (t_{j+1} - t_j)\|\nabla L(\phi_{z_2}(t_j,j))\|^2,
\end{align}
where the first equality uses~$\phi_{z_1}(t_{j+1},j+1) = \phi_{z_1}(t_{j+1},j) = \phi_{z_2}(t_j,j) - (t_{j+1} - t_j) \nabla L(\phi_{z_2}(t_j,j))$ from~\eqref{eq:hlast} and~\eqref{eq:z1update}, and the last equality uses~$\phi_{z_2}(t_{j+1},j) =\phi_{z_2}(t_j,j)$. Since~$t_{j+1}\!-\!t_j \leq \tau_{max} \leq \frac{\beta^2}{3K^3} \leq \frac{1}{3K} \leq \frac{2}{K}$, we have the desired property:~$L(\phi_{z_1}(t_{j+1},j+1)) \leq L(\phi_{z_2}(t_{j+1},j))$. Thus,
\begin{align}
    V(G(\phi(t_{j+1},j))) &- V(\phi(t_{j+1},j)) \leq 0. \label{eq:jinbound}
\end{align}

Following the work done in~\cite{chai17} and~\cite{goebel12}, we are able to perform direct integration in order to upper bound~$V(\phi(t,j))$ in terms of~$V(\phi(0,0))$ using~\eqref{eq:tinbound} and~\eqref{eq:jinbound} as bounds. Thus,
\begin{align}
    V(\phi(t,j)) \leq \exp \Big( -\frac{\beta A B}{2 K^2} t\Big) V(\phi(0,0)).
\end{align}
Using the comparison functions given in Lemma~\ref{lem:comparison}, we get a bound for~$|\phi(t,j)|^4$:
\begin{align}
    |\phi(t,j)|^4_\mathcal{A} &\leq \frac{16}{\beta^2} \exp \Big( -\frac{\beta A B}{2 K^2} t\Big) V(\phi(0,0)) \\
    &\leq \frac{8 K^2}{\beta^2} \exp \Big( -\frac{\beta A B}{2 K^2} t\Big) |\phi(0,0)|^4_\mathcal{A}.
\end{align}
Taking the fourth root gives the final answer. \hfill $\blacksquare$}{See~\cite{hendrickson21}. \hfill $\blacksquare$}

In practice, this preliminary result is useful when agreeing on initial values is easy to implement. However, it does not show our desired result, namely, exponential convergence to the minimizer of~$L$, regardless of initialization. By examining all possible scenarios at the first jump, we show in Theorem~\ref{thm:ges} below that exponential convergence to the minimizer of~$L$ still applies after the first jump.

\begin{theorem}[Exponential Convergence] \label{thm:ges}
Let Assumption~\ref{as:L} hold and consider the hybrid system~$\mathcal{H}$ defined in~\eqref{eq:hfirst}-\eqref{eq:hlast}. Let~$\mathcal{A}$ be as defined in~\eqref{eq:Adef} and choose~$\tau_{min}$ and~$\tau_{max}$ such that~$0 < \tau_{min} \leq \tau_{max} < \frac{\beta^2}{3K^3}$, where~$\beta$ is the strong convexity constant of~$L$ and~$K$ is the Lipschitz constant of~$\nabla L$. For each solution~$\phi$ and for all~$(t,j) \in \textnormal{dom } \phi$ such that~$j \geq 1$, the following is satisfied:
\begin{align}
    |\phi(t,j)|_\mathcal{A} &\leq \frac{8}{3} \sqrt[4]{2} \sqrt{\frac{K}{\beta}}   \exp \Big( -\frac{\beta A B}{8 K^2} t\Big) |\phi(0,0)|_\mathcal{A},
\end{align}
where~$A=\beta^2(1 - 2 \tau_{max} K) - \tau_{max} K^3 > 0$ and~$B=(1-2\tau_{max} K) \in (0,1)$. 
\end{theorem}
\noindent \emph{Proof: } \ifbool{Report}{ Two initialization scenarios must be considered: $\phi_{z_1}(0,0) = \phi_{z_2}(0,0)$ and~$\phi_{z_1}(0,0) \neq \phi_{z_2}(0,0)$. 

For the first case,~$\phi_{z_1}(0,0) = \phi_{z_2}(0,0)$, Proposition~\ref{prop:conv} applies in its original form. This is the best-case scenario that results in the smallest upper bound.

Now consider the second case, when~$\phi_{z_1}(0,0) \neq \phi_{z_2}(0,0)$. We note that after the first jump, all assumptions of Proposition~$\ref{prop:conv}$ hold. Thus, for any solution~$\phi$ where~$(t,j) \in \textnormal{dom } \phi$ such that~$j\geq 1$, we have the following: 
\begin{align}
    |\phi(t,j)|_\mathcal{A} &\leq \sqrt{\frac{K}{\beta}} \sqrt[4]{8} \exp \Big( -\frac{\beta A B}{8 K^2} t\Big) |\phi(t_1,1)|_\mathcal{A}, \label{eq:jgeq1}
\end{align}
where~$t_1$ denotes the time of the first jump.

We now seek to bound~$|\phi(t_1,1)|_\mathcal{A}$ in terms of~$|\phi(0,0)|_\mathcal{A}$. We begin by expanding and simplifying~$|\phi(t_1,1)|_\mathcal{A}^2$ using that~$\phi_{z_2}(t_1,1) = \phi_{z_1}(t_1,1)$ after the first jump:
\begin{align}
    |\phi(t_1,1)|_\mathcal{A}^2 &= \|\phi_{z_1}(t_1,1) - x^*\|^2 + \|\phi_{z_2}(t_1,1) - x^*\|^2= 2 \|\phi_{z_1}(t_1,1) - x^*\|^2, \label{eq:jumpnorm}
\end{align}
where~$x^*$ is the unique fixed point of~$\nabla L$.
Applying the flow dynamics and using~$\phi_{z_1}(t_1,1) = \phi_{z_1}(t_1,0)$ allows us to rewrite the distance of~$\phi_{z_1}$ from the minimizer of~$L$ after the first jump as
\begin{align}
    \|\phi_{z_1}(t_1,1) - x^*\|^2\!&=\!\|\phi_{z_1}(0,0)\!-\!t_1 \nabla L(\phi_{z_2}(0,0))\!-\!x^* \|^2 \\
    &\leq\|\phi_{z_1}(0,\!0)\!-\!x^* \|^2\!+\!2 t_1 \|\phi_{z_1}(0,\!0)\!-\!x^*\| \|\nabla L(\phi_{z_2}(0,\!0))\| + t_1^2 \|\nabla L(\phi_{z_2}(0,0)) \|^2 \\
    &\leq\|\phi_{z_1}(0,\!0)\!-\!x^* \|^2\!+\!2 t_1 K \|\phi_{z_1}(0,\!0)\!-\!x^*\| \|\phi_{z_2}(0,\!0)\!-\!x^*\| + t_1^2 K^2 \|\phi_{z_2}(0,0) - x^*\|^2, \label{eq:jumpexp}
\end{align}
where the last inequality comes from the~$K$-Lipschitz property of~$\nabla L$ and the fact that~$\nabla L(x^*) = 0$.
If~$\|\phi_{z_1}(0,0) - x^* \| \leq \|\phi_{z_2}(0,0) - x^* \|$, then~\eqref{eq:jumpexp} becomes
\begin{align}
    \|\phi_{z_1}(t_1,1)\!-\!x^*\|^2 &\leq (1\!+\!2t_1K\!+\!t_1^2K^2) \|\phi_{z_2}(0,0)\!-\!x^*\|^2 \\
    &\leq (1+t_1K)^2 |\phi (0,0)|_\mathcal{A}^2.
\end{align}
Similarly, if~$\|\phi_{z_1}(0,0) - x^* \| > \|\phi_{z_2}(0,0) - x^* \|$,  then~\eqref{eq:jumpexp} becomes~$\|\phi_{z_1}(t_1,1)\!-\!x^*\|^2 < (1+t_1K)^2 |\phi (0,0)|_\mathcal{A}^2$.
Thus, in both cases,
\begin{align}
    \|\phi_{z_1}(t_1,1)\!-\!x^*\|^2 \leq (1+t_1K)^2 |\phi (0,0)|_\mathcal{A}^2 \label{eq:firstflow}
\end{align}
Applying this to~\eqref{eq:jumpnorm} and using~$t_1 \leq \tau_{max} < \frac{\beta^2}{3K^3} < \frac{1}{3K}$, we get the bound
\begin{align}
    |\phi(t_1,1)|^2_\mathcal{A} &\leq 2 \left(\frac{4}{3}\right)^2 |\phi (0,0)|_\mathcal{A}^2. \label{eq:phit1bound}
\end{align}
Taking the square root and applying to~\eqref{eq:jgeq1} gives the final result for any~$(t,j) \in \textnormal{dom } \phi$ such that~$j\geq 1$. \hfill $\blacksquare$}{See~\cite{hendrickson21}. \hfill $\blacksquare$}

\section{Numerical Validation} \label{sec:sim}
We consider~$N=n$ agents for various values of~$n$. Each agent updates a scalar and they minimize
\begin{align}
    L(x) = \frac{1}{2}x^T Q x + b^Tx,
\end{align}
where~$x\in \R^n$,~$Q$ is a~$n \times n$ symmetric, positive definite matrix, and~$b\in \R^n$. 
To form Q, we decompose a random~$n \times n$ matrix into an unitary orthogonal matrix~$U$ and a matrix~$D$ that contains only our desired eigenvalues on the diagonal. We use these two resulting matrices to set~$Q = U^T D U$. Our choice of eigenvalues varies by trial (discussed below) with the minimum eigenvalue corresponding to~$\beta$ and the maximum eigenvalue corresponding to~$K$. The entries of~$b$ are
set to random values between 1 and 5. 
Simulations used 
the HyEq Toolbox (Version 2.04)~\cite{HyEqToolbox}\footnote{Simulation code for this section may be found at www.github.com/kathendrickson/DistrHybridGD.}. 

\begin{figure}[t!]
\centering
\includegraphics[width=8.8cm, trim={.25cm, .25cm, .25cm, .25cm},clip]{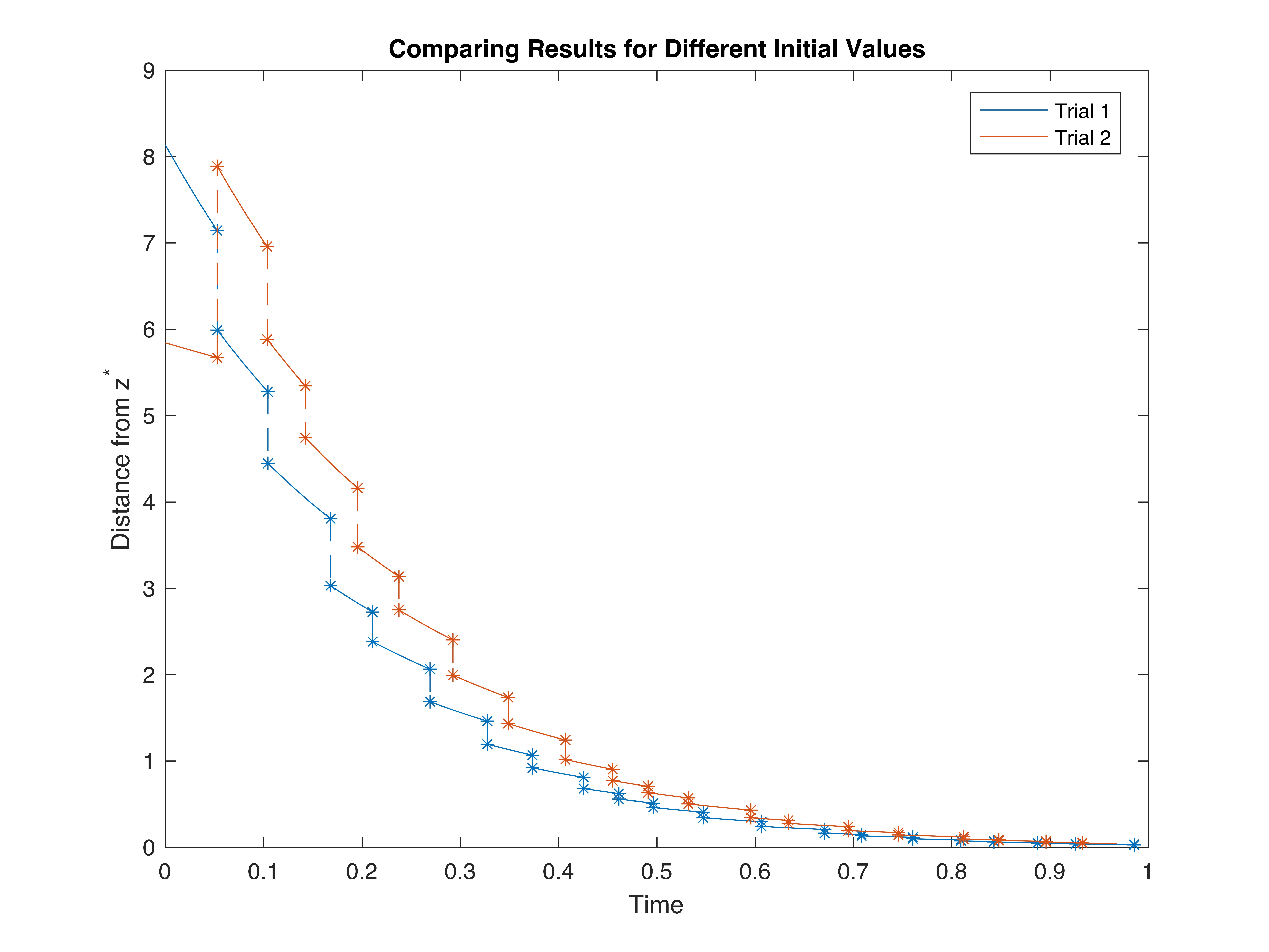}
\caption{Effect of initial values on convergence for two trials, where flows are denoted with solid lines and jumps with stars and dashed lines. Trial 1 sets~$\phi_{z_1} (0,0) = \phi_{z_2} (0,0)$, while Trial 2 sets~$\phi_{z_2} (0,0) = x^*$ instead. The first jump does not increase the distance from the minimizer in Trial 1 but does increase this distance in Trial 2. However, after the first jump, progress continues toward the optimum and differences between trials diminish.}
\label{fig:jump}
\end{figure}

We first compare convergence results for different initial values of~$\phi_{z_1}$ and~$\phi_{z_2}$ for five agents. For the first trial, we consider the case where~$\phi_{z_1} (0,0) = \phi_{z_2} (0,0) = (2,2,2,2,2)^T$. In Trial 2, we consider the ``worst-case'' initialization scenario: when~$\phi_{z_1} (0,0) = (2,2,2,2,2)^T$ is some distance from the optimum but~$\phi_{z_2} (0,0) = x^*$, resulting in an increase in the distance from the minimizer of~$L$ before the first jump. We consider~$\beta = K = 5$ and set~$\tau_{max} = \frac{\beta^2}{3K^3 + 1}$ and~$\tau_{min} = \frac{1}{2}\tau_{max}$. Figure~\ref{fig:jump} shows the distance from optimum through the first twenty jumps for both trials. There is a consistent decrease in the distance to the minimizer, even at jumps, for the first trial. In contrast, when initial values for~$\phi_{z_1}$ and~$\phi_{z_2}$ are not equal, there is an increase in distance to the minimizer after the first jump in the second trial. However, as expected, distance to the optimum
decreases exponentially thereafter, with the difference between the two trials decreasing over time.

\begin{figure}[t!]
\centering
\includegraphics[width=8.8cm, trim={.25cm, .25cm, .25cm, .25cm},clip]{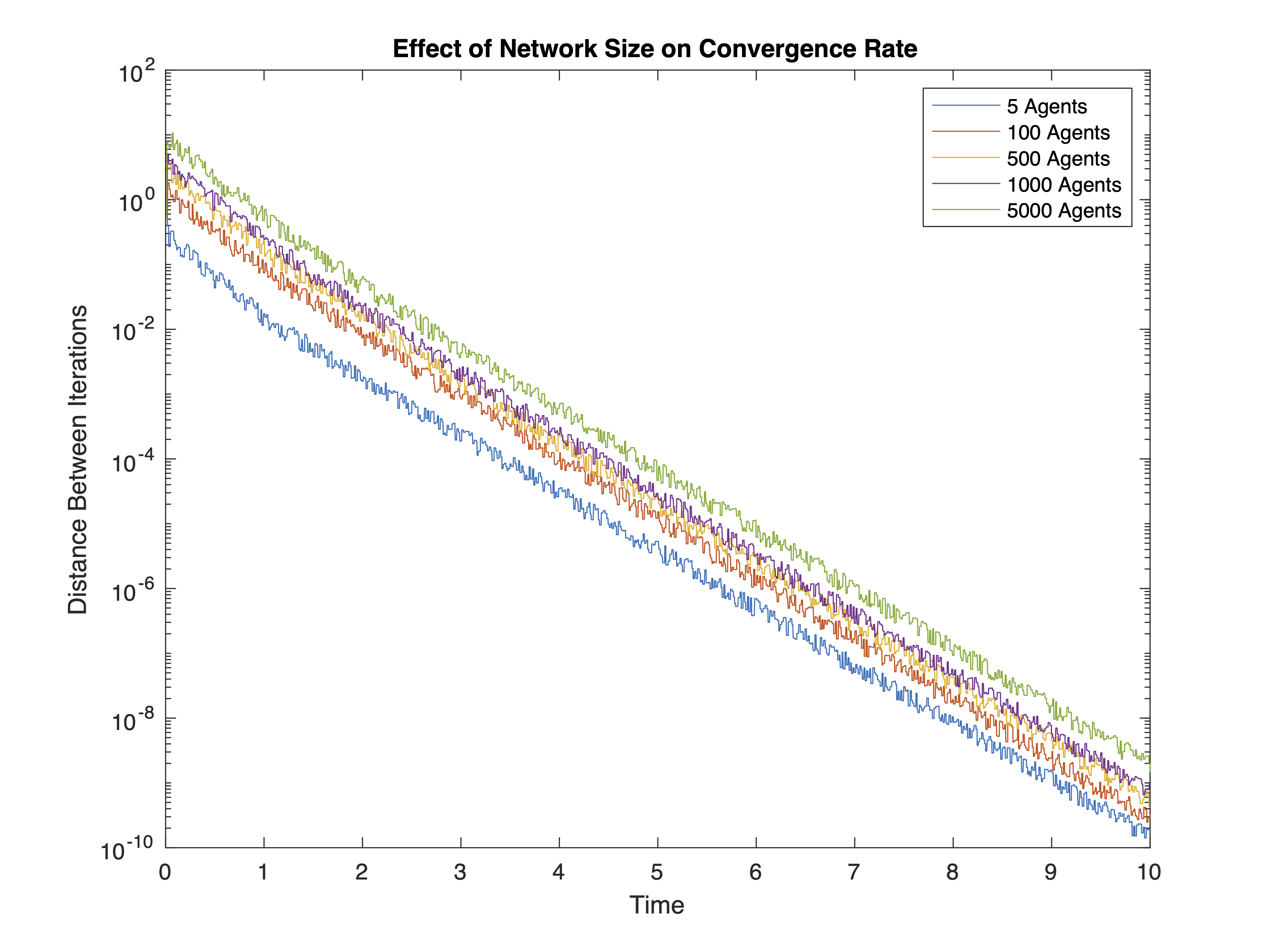}
\caption{Effect of network size on convergence. Convergence results still hold even for very large network sizes, demonstrating the scalability of our algorithm.}
\label{fig:networksize}
\end{figure}
\balance
We then examined the effects of varying the network size from 5 agents to 100, 500, 1000, and 5000 agents. We set~$\beta = 2$ and~$K = 4$ and chose to initialize~$\phi_{z_1}$ and~$\phi_{z_2}$ with vectors of twos in~$\R^n$. For each network size, the matrix~$Q$ and vector~$b$ were randomly generated. As shown in Figure~\ref{fig:networksize}, drastically expanding the network size does not have a significant impact on convergence. This demonstrates our algorithm's scalability and convergence results that hold regardless of network size. 
 
\section{Conclusion} \label{sec:concl}
This paper presented a hybrid systems framework for analyzing continuous-time multi-agent optimization with discrete-time communications. Using this framework, we established that every maximal solution is complete, as well as the exponential convergence of a block coordinate descent law to the minimizer of a strongly convex and smooth objective function. Future work in this area includes the use of heterogeneous timers and exploration of other update laws, as well as constrained problems.

\ifbool{Report}{\appendix
\emph{Proof of Lemma~\ref{lem:solutions}:} Using Proposition~6.10 in~\cite{goebel12} with~$U=C$, let~$\xi=(z_1, z_2, \tau) \in C \backslash D$. Then~$f(\xi) \subset T_C(\xi)$. Because~$G(D) \subset C$, case (c) in Proposition~6.10 does not apply. We avoid case (b) of Proposition~6.10 by showing that every solution lies entirely in a compact subset~$W \subset C$. To do this, we verify is that there is no finite escape time for any solution. Consider a solution~$\phi$. Then~$\phi_{z_1}(0,0)$ and~$\phi_{z_2}(0,0)$ denote the initial values of~$\phi_{z_1}$ and~$\phi_{z_2}$, respectively. Let~$t_1$ denote the continuous-time at which the first jump occurs; then the value of~$\phi$ after the first jump may be written as~$\phi(t_1,1)$. We first show that there is no finite escape time from initialization through the first jump. Towards doing this, we use~\eqref{eq:phit1bound},~$|\phi(t_1,1)|^2_\mathcal{A} \leq 2 (\frac{4}{3})^2 |\phi (0,0)|_\mathcal{A}^2$, and apply the comparison functions in Lemma~\ref{lem:comparison}. This gives the set of inequalities
\begin{align}
    V(\phi(t_1,1)) &\leq \frac{K^2}{2}|\phi(t_1,1)|^4_\mathcal{A} \\
    &\leq \frac{K^2}{2} 4 \left(\frac{4}{3}\right)^4 |\phi (0,0)|_\mathcal{A}^4 \\
    &\leq \frac{K^2}{2} 4 \left(\frac{4}{3}\right)^4 \frac{16}{\beta^2} V(\phi(0,0)),
\end{align}
which follow from~\eqref{eq:compupper},~\eqref{eq:phit1bound}, and~\eqref{eq:complower}, respectively.
Thus, through the first jump,~$V(\phi(t_1,1)) \leq \frac{8192 K^2}{81\beta^2} V(\phi(0,0))$.

After this first jump,~$\phi_{z_1}(t_1,1) = \phi_{z_2}(t_1,1)$ then holds. By construction and Assumption~\ref{as:L},~$\nabla L$ is Lipschitz and thus the map~$f$ is Lipschitz as well. Applying~\eqref{eq:dotVupper} after hybrid time~$(t_1,1)$, we have
\begin{align}
    \dot{V}(\phi(t,j)) = \langle \nabla V(\phi(t,j)), f(\phi(t,j)) \rangle \leq 0, \label{eq:dotvneg}
\end{align}
for all~$(t,j) \in \textnormal{dom }\phi$ such that~$j \geq 1$.
Thus, for any solution~$\phi$, we see that
\begin{equation} \label{eq:vjumpneg}
V(\phi(t,j)) \leq \frac{8192 K^2}{81 \beta^2} V(\phi(0,0)) 
\end{equation}
for all~$(t,j) \in \textnormal{dom }\phi$. 

Now consider a solution to~$\dot{\xi} = f(\xi)$ that starts from some~$c$-sublevel set~$W_1 = \{\xi \in \mathcal{X} : V(\xi) \leq c\}$. Then, from~\eqref{eq:dotvneg} and~\eqref{eq:vjumpneg}, we see that 
all such solutions remain in the sub-level set~$W_2 = \{\xi \in \mathcal{X} : V(\xi) \leq \max\{c, \frac{8192K^2}{81\beta^2}V\big(\phi(0, 0\big)\}\}$. 
Because~$V$ is continuous and radially unbounded by Lemma~\ref{lem:comparison},~$W_2$ is compact. From~\eqref{eq:dotvneg},~$W_2$ is forward invariant for~$\mathcal{H}$. Thus any trajectory that starts in the subset~$W_1$ remains in~$W_2$. Thus by Theorem 3.3 in~\cite{khalilnonlinear}, there is no finite escape time from~$C$.
\hfill $\blacksquare$}{}

\bibliographystyle{IEEEtran}
\bibliography{IEEEabrv,sources}

\end{document}